\documentclass[11pt,a4paper]{elsarticle}
\usepackage{amssymb,theorem,amsmath}
\newtheorem{thm}{Theorem}[section]
\newtheorem{lem}{Lemma}[section]
\newproof{pf}{\bf Proof}

{\theorembodyfont{\normalfont}
\newtheorem{rmk}{Remark}[section]}
\newtheorem{cor}{Corollary}[section]
\newdefinition{defn}{Definition}[section]
{\theorembodyfont{\normalfont}
\newtheorem{exmp}{Example}[section]}

\makeatletter

\let \al=\alpha
\let \be=\beta

\let \vare=\varepsilon

\let \de=\delta
\let \th=\theta
\let \la=\lambda

\let \ga=\gamma
\let \p=\partial
\let \q=\quad

\let \med=\medskip
\let \smal=\smallskip
\let \dps=\displaystyle

\let \ol=\overline

\newcommand{\R}{\mathbb{R}}
\newcommand{\N}{\mathbb{N}}

\def\system#1{\left\{\null\,\vcenter{
\ialign{\strut\hfil$##$&$##$\hfil&&\enspace$##$\enspace&
\hfil$##$&$##$\hfil\crcr#1\crcr}}\right.}

\begin{document}

\begin{frontmatter}

\title{Periodic solutions  for a non-monotone family of 
delayed differential equations with applications to Nicholson systems }
\author{Teresa Faria{\footnote{Tel:~+351 217500192,
fax:~+351 217500072, e-mail:~teresa.faria@fc.ul.pt.}}}

\address{Departamento de Matem\'atica and CMAF-CIO, Faculdade de Ci\^encias, Universidade de Lisboa\\ Campo Grande, 1749-016 Lisboa, Portugal}



\begin{abstract} For a  family of $n$-dimensional periodic delay differential equations which encompasses a broad set of models used  in structured population dynamics,  the existence of a positive periodic solution  is obtained under very mild conditions. The proof uses the Schauder fixed point theorem and  relies on the  permanence of the system.
A  general criterion for the existence of a positive periodic solution for  Nicholson's blowflies periodic systems (with both distributed and discrete time-varying delays) is derived as a simple application of our main result, generalizing the few existing results concerning multi-dimensional  Nicholson  models. 
  In the case of a Nicholson system with  discrete delays all  multiples of the period,  the global attractivity of the positive periodic solution is further analyzed, improving results in recent literature.  
 \end{abstract}


\end{frontmatter}

\section{Introduction}

 In recent years, the question of the existence of periodic solutions for  periodic delay differential equations (DDEs)  has attracted the interest of many researchers, and a plethora of  positive answers has been provided by using a variety of methods. To a large extent, the techniques used in the literature apply  to a specific equation  only, while other ones apply to a  very particular class of   DDEs, with emphasis  on scalar models.
 For some classical  models from mathematical biology, the available existence results require a very restrictive set of assumptions, which are not easily verifiable, much less  extendable to other families of DDEs.

The main purpose of this paper is to investigate  the existence of a positive periodic solution for a broad  class of periodic and in general non-monotone $n$-dimensional DDEs which encompasses a  large number of population models with patch structure. DDEs with patch structure 
have extensive applications  in  population dynamics, where the patch structure accounts for  situations of  heterogeneous environments due to several aspects,  or in disease and epidemic models with  different classes for cells or individuals, with transition among the classes. In particular, the study of periodic models is especially significant,  as they reflect    periodical variations of the weather or  seasonality of the habitat in general, so the quest for positive periodic solutions for such models becomes quite relevant.

In this paper, we  consider a family of periodic delayed population models with patch structure and multiple time-varying delays of the form 
\begin{equation}\label{1.0}
x_i'(t)=-d_i(t)x_i(t)+\sum_{j=1,j\ne i}^n a_{ij}(t)x_j(t)+\sum_{k=1}^m \be_{ik}(t)   \int_{t-\tau_{ik}(t)}^t\!\! b_{ik}(s,x_i(s))\, d_s\eta_{ik}(t,s),\  i=1,\dots,n,
\end{equation}
where all the  coefficients and delay functions are assumed to be continuous,   non-negative and periodic on $t$, with a common period $\omega>0$, and $\eta_{ik}(t,s)$ are  bounded, nondecreasing  on $s$, locally integrable and $\omega$-periodic on $t$. Some additional conditions on the  coefficients $d_i(t), a_{ij}(t),\be_{ik}(t)$ and on the nonlinearities $b_{ik}(t,x)$ will be assumed.  Special attention will be given to the study of \eqref{1.0} with  $\eta_{ik}(t,s)=H_{{t-\tau_{ik}(t)}}(s)$, where $H_t(s)$ is the Heaviside function $H_t(s)=0$ if $s\le t$, $H_t(s)=1$ if $s> t$. In this case, we obtain a system with discrete delays of the form
\begin{equation}\label{1.1}
x_i'(t)=-d_i(t)x_i(t)+\sum_{j=1,j\ne i}^n a_{ij}(t)x_j(t)+\sum_{k=1}^m \be_{ik}(t)  h_{ik}(t,x_i(t-\tau_{ik}(t))),\  i=1,\dots,n.
\end{equation}

Many important delayed non-autonomous models from mathematical biology can be written in the form  \eqref{1.0}, see e.g.~\cite{Kuang,MetzDiek,Smith11}. In Section 2, a descriptive set of hypotheses, as well as  a brief biological interpretation of the model, will be given.

 The present paper is a continuation of the  research recently conducted by Faria, Obaya and Sanz in \cite{FOS}, where  the asymptotic behavior of  solutions
     for non-autonomous systems   \eqref{1.1} 
  was carefully analyzed, and sufficient conditions for either the permanence or extinction of all population  given. A permanence result was established in \cite{FOS} for a generic  system \eqref{1.1} with all the  coefficients and delays  given by non-negative, continuous, bounded  functions (not necessarily periodic), under very mild and optimal sufficient conditions. As we shall see, the permanence result in \cite{FOS} can be easily  extended to systems \eqref{1.0}. Here, the leading  ideas are, on one hand,  to interpret \eqref{1.0} as the result of adding bounded delayed perturbations to a linear homogeneous ordinary differential equation (ODE)
 and, on the other hand,  to  use the uniform persistence of   \eqref{1.0}:
  under the same  hypotheses for permanence, by applying  the Schauder fixed point theorem we  further show that at least one positive period solution must exist.
 Our results are achieved under a very weak set of assumptions and have significant applications.

 Among them, and as an important illustration, we have in mind 
to apply our results to   Nicholson systems
with  time-dependent, either distributed or discrete, delays, of the forms
\begin{equation}\label{1.2'}
x_i'(t)=-d_i(t)x_i(t)+\sum_{j=1,j\ne i}^n a_{ij}(t)x_j(t)+\sum_{k=1}^m \be_{ik}(t)  \int_{t-\tau_{ik}(t)}^t\!\!  \ga_{ik}(s)x_i(s)e^{-c_{ik}(s)x_i(s)}\, ds, \  i=1,\dots,n,
\end{equation}
or
\begin{equation}\label{1.2}
x_i'(t)=-d_i(t)x_i(t)+\sum_{j=1,j\ne i}^n a_{ij}(t)x_j(t)+\sum_{k=1}^m \be_{ik}(t)    x_i(t-\tau_{ik}(t))e^{-c_{ik}(t)x_i(t-\tau_{ik}(t))}, \  i=1,\dots,n,
\end{equation}
where  the coefficients and delays  are continuous and bounded on $\R$, with $d_i(t),c_{ik}(t)>0, a_{ij}(t), \be_{ik}(t),\ga_{ik}(t), \tau_{ik}(t)\ge 0 $. 
In view of our purposes, we give here some references on Nicholson's equations and systems, with emphasis on periodic versions of such models.


The {\it Nicholson' blowflies equation}
\begin{equation}\label{N}
N'(t)=-dN(t)+\be N(t-\tau)e^{-aN(t-\tau)}\q (d,\be, a,\tau>0)
\end{equation}
was introduced by Gurney et al.~in 1980 \cite{GBN}, and its  biological  impact  was immediately apparent, as the proposed model agreed with  Nicholson's experimental data  on the Australian sheep blowfly (see e.g.~\cite{Nicholson}). 
Since then, an immense literature concerning Nicholson's equation,    generalizations,   related models and applications to real world problems has been produced. 
The periodic version of  \eqref{N}, given by
\begin{equation}\label{NP}
N'(t)=-d(t)N(t)+\be(t) N(t-\tau(t))e^{-a(t)N(t-\tau(t))},
\end{equation}
 with $d(t),\be(t),\tau(t),a(t)$ positive, $\omega$-periodic continuous functions ($\omega>0$),   has been studied in a number of papers. 
  The simplest version of \eqref{NP} is
 when the delay is a multiple of the period,  in which case the $\omega$-periodic solutions of \eqref{NP} are exactly the $\omega$-periodic solutions of the equation with no time delay $N'(t)=-d(t)N(t)+\be(t) N(t)e^{-a(t)N(t)}$. For the particular situation of $\tau(t)=m\omega\ (m\in\N)$ and $a(t)\equiv a>0$,   Saker and Agarwal \cite{SakerAg} showed that there is a  positive $\omega$-periodic solution $N^*(t)$ of \eqref{NP} if $
 \min_{t\in [0,\omega]} \be(t)>\max_{t\in [0,\omega]} d(t),$
and  gave some additional  conditions for its global  attractivity. See also \cite{LiDu} for a refinement of the result in \cite{SakerAg}.
A significant breakthrough was later  achieved by Chen \cite{Chen}, who used the continuation theorem of coincidence degree  to establish the existence  of a positive $\omega$-periodic solution of \eqref{NP} under much more general conditions.
More recently,  an elegant unifying method, based on  the continuation theorem, was proposed by Amster and Idels \cite{AmsterIdels}  to show the existence of positive periodic solutions for a general class of scalar period DDEs with the form $x'(t)=\mp a(t)f(x(t))x(t)\pm \la b(t)g(x_t)$ ($\la>0$ a parameter). Their results apply to the  case of a Nicholson scalar equation with distributed delay, as well as
 to other  important biological models.  For other   criteria of  existence of periodic solutions for scalar periodic DDEs  within the class $x'(t)=\mp a(t)f(x(t))x(t)\pm \la b(t)g(x_t)$ and based on several fixed point methods,  see  \cite{DingNieto,FLT,Liu14,Wang04,ZDC} and references therein. 

 Only recently has some attention  been given to multi-dimensional  versions of Nicholson models, with priority in  {\it autonomous} Nicholson systems  \cite{BIT,Faria11,FariaRost,Liu}. 
 For $n>1$,
 very little is known about positive periodic solutions for general
 non-monotone periodic DDEs. In what concerns  periodic Nicholson systems, results concerning the existence of positive periodic solutions  in the case $n=2$ have been established in a few papers, see \cite{Liu11,Troib}, and seem to be virtually non-existent for the situation of  \eqref{1.2} with $n>2$, or for the case of distributed delays  \eqref{1.2'} with $n\ge 2$.
 For  related results  for periodic or almost-periodic Nicholson systems with harvesting terms, we refer also  to \cite{ Wang13,WWC,Zhou}.


  In spite of the variety of methods and tools that have been proposed, to the best of our knowledge,  there is no general result in the literature concerning the existence of positive periodic solutions for periodic $n$-dimensional Nicholson systems  \eqref{1.2'} or  \eqref{1.2}. 
 Surprisingly, here a general criterion is obtained as an immediate consequence of our  Theorem \ref{thm3.1},  established for a far more general framework.  
 The results in \cite{LiDu}  and  \cite{Liu14} are recovered by our general criterion, when applied to the   scalar version of \eqref{1.2}.   We however believe that sharper results are to be expected, under more natural restrictions involving the average integrals of the coefficients over the interval $[0,\omega]$, as in \cite{Chen, Troib} for $n=1$ -- rather than the  pointwise values of such coefficients, as in the results presented here. This will be the subject of future research. 

The contents  of the remainder of the paper are  now summarized.   Section 2 is a section of preliminaries,  where a set of assumptions  is introduced and the general criterion for permanence in \cite{FOS}  extended to the family of  DDEs \eqref{1.0}. The main result of the paper, Theorem \ref{thm3.1}, is given in Section 3:  in the case of {\it periodic} systems \eqref{1.0} we show that the
 sufficient conditions for permanence  are enough to guarantee the existence of at least one positive $\omega$-periodic  solution. The result for  periodic $n$-dimensional Nicholson systems \eqref{1.2'}  is deduced as a particular case.  In Section 4, we consider  \eqref{1.2}  with all the delays multiple of the period, and give sufficient conditions for the global attractivity of the positive periodic solution. Our results extend the ones in \cite{LiDu, Liu14} and improve some criteria in \cite{AmsterIdels,Troib}. The situation of systems with autonomous coefficients is also considered and  the global asymptotic stability of a positive equilibrium deduced under optimal conditions, generalizing the results in \cite{BIT,FariaRost}.  
 Although emphasis  is given  to periodic Nicholson systems,  other relevant population models satisfy the  hypotheses imposed here; see Sections 2 and 3 for  examples.
 
 \section{Preliminaries }
\setcounter{equation}{0}

We start by introducing some standard notation. For $\tau\ge 0$, set $C:=C([-\tau,0];\R^n)$ to be the Banach
space   endowed with the norm
$\|\phi\|=\max_{\th\in[-\tau,0]}|\phi(\th)|$, where $|\cdot|$ is  a fixed norm in $\R^n$.  We shall also use $|A|$ to denote the (operator) norm of an $n\times n$ matrix $A$ with constant entries.  A vector $v\in \R^n$  is identified in $C$ with the constant function $\psi(s)=v$ for $-\tau\le s\le 0$. 

A DDE in $C$ takes the general form 
\begin{equation}\label {DDE}
x'(t)=f(t,x_t),
\end{equation}
 where $f:\Omega\subset \R\times C\to \R^n$ and $x_t$ denotes the restriction of a solution $x(t)$ to the time interval $[t-\tau,t]$, i.e.,  $x_t\in C$ is given by $x_t(\th)=x(t+\th), -\tau\le \th\le 0$. Take $\Omega=[\al,\infty)\times D$ with $\al\in\R$  and
 $D\subset C$,  and suppose that
 $f$ is continuous and regular enough  so that 
the initial value problem is well-posed, in the sense that for each $(\sigma,\phi)\in [\al,\infty)\times D$ there exists a unique solution of  the problem $x'(t)=f(t,x_t), x_\sigma=\phi$,  defined on a maximal interval of existence: in this situation, this solution is denoted by $x(t,\sigma,\phi)$ in $\R^n$ or  $x_t(\sigma,\phi)$ in $C$.   Whenever necessary, the more explanatory notation  $x(t,\sigma,\phi,f)$ is used.

 We designate by $C^+$  the cone of nonnegative functions in $C$, $C^+=C([-\tau,0];[0,\infty)^n)$, and by  $int \, C^+$ its interior. In $C$, $\le $ denotes the usual partial order generated by  $C^+$: $\phi\le \psi$ if and only if $\psi-\phi\in C^+$; by $ \phi\ll \psi $,   we mean that $\psi-\phi\in int\, C^+$. The relations $\ge$ and $\gg $ are defined in the obvious way; thus,  we write $\psi\ge 0$ for $\psi\in C^+$ and $\psi\gg 0$ for $\psi \in int \, C^+$.  
The situation of no delays ($\tau=0$) is included in our setting, in which case $C$ is identified with $\R^n$ and
 $C^+$ with $\R_+^n:=[0,\infty)^n$. 

For simplicity, here we say that  \eqref{DDE} (or $f$) is {\it cooperative} if it satisfies
Smith's {\it quasimonotone condition},  given by (see \cite{Smith})

\vskip 2mm

{(Q)} for $\phi,\psi\in D,\phi\le \psi$ and $\phi_i(0)=\psi_i(0)$, then $f_i(t,\phi)\le f_i(t,\psi), \ i=1,\dots,n,t\ge \al$.

\vskip 2mm

\noindent Condition  (Q)   allows   comparison of solutions between two related DDEs $x'(t)=f(t,x_t)$ and $x'(t)=g(t,x_t)$:  if $f\le g$ on $[\al,\infty)\times D$ and either $f$ or $g$ is cooperative, then, for $\sigma\ge \al$ and $\phi, \psi\in D$ with $\phi\le \psi$, we have
$x(t,\sigma, \phi,f)\le x(t,\sigma, \psi,g)$ for $t\ge \sigma$ whenever the solutions are defined (see \cite{Smith}).
In particular, (Q) guarantees the monotonocity of solutions of   \eqref{DDE} relative to initial data.

\med

Consider a family of non-autonomous  systems \eqref{1.0}, 
  and further suppose that  $b_{ik}(t,0)=0$ for $t\in\R$  and $b_{ik}$ have
 partial derivatives with respect to the second variable at $x=0^+$,  given by 
$ \frac{\p b_{ik}}{\p x}(t,0)=\ga_{ik}(t)$, for all $i,k$.
  In this way, system \eqref{1.0} takes the general form
  \begin{equation}\label{0}
  \begin{split}
x_i'(t)=-d_i(t)x_i(t)&+\sum_{j=1,j\ne i}^n a_{ij}(t)x_j(t)\\
&+\sum_{k=1}^m \be_{ik}(t)  \int_{t-\tau_{ik}(t)}^t\ga_{ik}(s)h_{ik}(s,x_i(s))\,  d_s\eta_{ik}(t,s),\  i=1,\dots,n,
\end{split}
\end{equation}
 where all the coefficients, kernels and delay functions  are supposed to be continuous, bounded and nonnegative. 
 As a special case of \eqref{0}, we shall consider systems with  time-dependent discrete delays in the nonlinear terms,  written in the form
\begin{equation}\label{3}
x_i'(t)=-d_i(t)x_i(t)+\sum_{j=1,j\ne i}^n a_{ij}(t)x_j(t)+\sum_{k=1}^m \be_{ik}(t)    h_{ik}(t,x_i(t-\tau_{ik}(t))),\  i=1,\dots,n.
\end{equation}

Systems \eqref{0} and  \eqref{3} are considered as abstract DDEs  in the phase
space $C=C([-\tau,0];\R^n)$, where
$$\tau=\sup \{\tau_{ik}(t): t\ge 0, i=1,\dots, n,\, k=1,\dots,m\}.$$  For future reference, $\R^n$ is supposed to be equipped with the supremum norm $|x|=\max_{1\le i\le n}|x_i|$, $x=(x_1,\dots,x_n)\in\R^n$.
Note that \eqref{0} is obtained by adding a delayed perturbation $M(t,x_{t})$   to the  linear   ODE
 \begin{equation}\label{2.1}
x_i'(t)=-d_i(t)x_i(t)+\sum_{j=1,j\ne i}^n a_{ij}(t)x_j(t),\q i=1,\dots,n,
\end{equation}
with  $M(t,x_{t})$ of the form
$M(t,x_{t})=(M_1(t,x_{1,t}),\dots, M_n(t,x_{n,t}))$; for \eqref{0}  each component $M_i(t,\phi_i)$ is given by
 \begin{equation}\label{M0}
 M_i(t,\phi_i)= \sum_{k=1}^m \be_{ik}(t) \int_{t-\tau_{ik}(t)}^t\ga_{ik}(s)h_{ik}(s,\phi_i(s-t))\,  d_s\eta_{ik}(t,s),\ 
 \end{equation}
 whereas for \eqref{3} the components $M_i(t,\phi_i)$  read as
 \begin{equation}\label{M1}
 M_i(t,\phi_i)= \sum_{k=1}^m \be_{ik}(t) h_{ik}(t,x(t-\tau_{ik}(t))),\ 
 \end{equation}
for $\ t\in\R, \phi_i\in C([-\tau,0];\R), i=1,\dots, n.$

Typically, system \eqref{0}  can be used to model the  population growth of either a single or multiple species  structured into $n$ classes or patches, with migration among them: $x_i(t)$ denotes the density of the $i$th-species population, $a_{ij}(t)$ is the dispersal rate of the population migrating  from class $j$ to class $i$, $d_i(t)$ is the coefficient of instantaneous loss for class $i$ (which integrates both the death rate and the migration coefficients referring to the individuals that leave class $i$ to move to other classes), and  $M_i(t,\phi_i)$  is the birth function for class $i$.
  Following the general approach in the literature -- though not always justifiable from a biological viewpoint \cite{Diekmann01} --  multiple time-varying  delays have been incorporated in  the birth contribution.  
 

\med
Throughout the paper, hypotheses  will be taken from the following set of conditions:

\begin{itemize}

\item[(H0)] the  functions $d_i,a_{ij},\be_{ik},\ga_{ik}, h_{ik}(\cdot, x) \, (x\ge 0), \eta_{ik}(\cdot, s)\, (s\in\R)$ and $\tau_{ik}$ are  $\omega$-periodic ($\omega>0$) on $t\in \R$; 

\item[(H1)]  $d_i,a_{ij}:\R\to\R\, (j\ne i)$ are   continuous, with $a_{ij}(t)\ge 0,i\ne j, d_i(t)>0$ for $t\in \R$ and $i,j\in\{1,\dots,n\}$;

\item[(H2)] there exist a vector $u=(u_1,\dots,u_n)\gg 0$ and $t_0\in \R$ such that
$d_i(t)u_i\ge \sum_{j=1,j\ne i}^n a_{ij}(t)u_j$ for $t\in\R$, with $d_i(t_0)u_i> \sum_{j=1,j\ne i}^n a_{ij}(t_0)u_j$, $i\in\{1,\dots,n\}$;

\item[(H3)]  $\tau_{ik},\be_{ik},\ga_{ik}:\R\to [0,\infty)$ are continuous, $\eta_{ik}:\R\times\R\to \R$ are  bounded, with $\eta_{ik}(t,s)$  nondecreasing on $s$ and locally  integrable on $t$, and
\begin{equation}\label{beta-i}
\be_i(t):= \sum_{k=1}^m \be_{ik}(t) \int_{t-\tau_{ik}(t)}^t\ga_{ik}(s)\, d_s\eta_{ik}(t,s)>0,\q {\rm}\q t\in \R,
\end{equation}
for $i\in \{1,\dots,n\},k\in \{1,\dots,m\};$

\item[(H4)] $h_{ik}:\R\times [0,\infty)\to [0,\infty)$ are bounded, continuous and  locally Lipschitzian in $x$, with $h_{ik}(t,0)=0$ for $t\in\R$ and
$$h_{ik}(t,x)\ge h_i^-(x)\q {\rm}\q t\in\R,x\ge 0,k=1,\dots,m,$$ where $h_i^-:[0,\infty)\to [0,\infty)$ is  continuous on $[0,\infty)$, continuously differentiable in a right  neighborhood of $0$,  with $h_i^-(0)=0,(h_i^-)'(0)=1$ and $h_i^-(x)>0$ for $x>0$,  $i\in \{1,\dots,n\}$.
 

\end{itemize}

Assumptions (H1)-(H4), together with either (H0) or the boundedness of all functions in their domains,  guarantee the existence and uniqueness of  solutions for the initial value problems of \eqref{0} with $x_\sigma=\phi\in C^+$, defined for $t\ge \sigma$ \cite{HaleLunel}. For \eqref{3}, $\be_i(t)$ in \eqref{beta-i} reads simply as 
$\be_i(t)=\sum_{k=1}^m \be_{ik}(t), i=1,\dots,n$. Typically, in (H4) we take $h_i^-(x)=\min\{ h_{ik}(t,x):  t\in [0,\omega], 1\le k\le m\}$.

Hereafter, we designate by $A(t),B(t),D(t),M(t)$  the $\omega$-periodic  $n\times n$ matrices defined on $\R$ and given by
\begin{equation}\label{2.3}
\begin{split}
D(t)=diag\, (d_1(t),\dots,d_n(t)),&\q  A(t)=[a_{ij}(t)]\\
 B(t)=diag\, (\be_1(t),\dots,\be_n(t)),&\q
 M(t)=B(t)+A(t)-D(t),
\end{split}
\end{equation}
where $a_{ii}(t)\equiv0,\, 1\le i\le n$. In the literature, $M(t)$ is often called the {\it community matrix} for \eqref{0}. In biological terms,  (H1)-(H3) are quite natural conditions for periodic structured population models; for a discussion see   \cite{FOS}, also for further references. Jointly with (H0)-(H4), we shall also consider the following assumption:
  \begin{itemize}
   \item[(H5)]   there exists $v=(v_1,\dots,v_n)\gg 0$ such that $M(t)v\gg 0$ for $t\in [0,\omega].$
\end{itemize}

Besides \eqref{1.2'}, where the  nonlinearities are of  Ricker-type, other useful population models  satisfying the above hypothesis (H4) can be considered. Among them, models  \eqref{0} with 
$h_{ik}(t,x)=xe^{-c_{ik}(t)x^\al}\, (\al>0)$
or  with nonlinearities of Mackey-Glass type
\begin{equation}\label{hiMG}
h_{ik}(t,x)=\frac{x}{1+c_{ik}(t)x^\al}\q (\al\ge1),
\end{equation}
where $c_{ik}(t)$ are continuous, positive and bounded, satisfy (H4). See Section 3 for an illustrative example.

 \smal

Motivated by its biological interpretation, 
only nonnegative solutions of \eqref{0} are meaningful, and therefore admissible. Here, initial conditions are taken in  $C_0$, where $$C_0=\{ \phi\in C^+: \phi(0)\gg 0\}.$$ 

The notions of  uniform persistence  and permanence given below (see e.g. \cite{Kuang})  will always refer to the choice of $C_0$ as the set of admissible initial conditions, given, as  convention,  at the instant of time $t=0$: i.e., initial conditions read as $x_0=\phi\in C_0$; of course, one can replace  $[0,\infty)$ by any  time interval $[\al,\infty),\al \in \R.$

\begin{defn}\label{def2.1} A DDE $x'(t)=f(t,x_t)$  
 is said to be {\it uniformly persistent} (in $C_0$)  if all solutions $x(t,0,\phi)$ with $\phi\in C_0$ are defined  on $[0,\infty)$ and there is  $m>0$ such that $\dps\liminf_{t\to\infty}x_i(t,0,\phi)\ge m$ for all $1\le i\le n,\phi\in C_0$.
The DDE $x'(t)=f(t,x_t)$   is said to be {\it permanent} (in $C_0$)  if it is dissipative and uniformly persistent; in other words, all solutions $x(t,0,\phi),\phi\in C_0$, are defined on $[0,\infty)$ and there are positive constants $m,L$ such that, given any $\phi\in C_0$, there exists $t_0=t_0(\phi)$  for which
\begin{equation}\label{3.5}
m\le x_i(t,0,\phi)\le L\q {\rm for}\q t\ge t_0,\ i=1,\dots,n.
\end{equation} 
\end{defn}

In general, the nonlinearities in \eqref{0} are non-monotone in $x$,  thus monotone techniques  do not apply directly. Nevertheless, for the case of systems with discrete delays \eqref{3},  in \cite{FOS} the authors considered convenient auxiliary cooperative systems, and exploited results from the theory of monotone DDEs  as in \cite{Smith}, to deduce the global asymptotic behavior of solutions. For {\it periodic} systems \eqref{0},
some   consequences and generalizations of results in \cite{FOS} are given in the following  theorem:

\begin{thm}\label{thm2.1} 
(i) If (H0)-(H2) are satisfied, then the $\omega$-periodic linear homogeneous system
\eqref{2.1}
 is cooperative and exponentially asymptotically stable.\vskip 0cm
  (ii)  If (H0)-(H4) are satisfied, all solutions of \eqref{0} with initial conditions in $C_0$ are defined and  strictly positive on $[0,\infty)$; moreover,  \eqref{0} is dissipative (in $C_0$). \vskip 0cm
   (iii)  If (H0)-(H5) are satisfied, then \eqref{0}
is permanent (in $C_0$). 
\end{thm}

\begin{pf} The assertions in (i) and (ii)  are immediate consequences of Theorems 2.1 and 2.3 in \cite{FOS}. For the proof of (iii), below we adapt the arguments for the proof of  \cite[Theorem 3.3]{FOS}, omitting however details.

After the scaling of variables   $\hat x_j(t)=x_j(t)/v_j$, system \eqref{0} reads as 
  \begin{equation*}
  \begin{split}
\hat x_i'(t)=- d_i(t)\hat x_i(t)&+\sum_{j=1,j\ne i}^n \hat a_{ij}(t)\hat x_j(t)\\
&+\sum_{k=1}^m \be_{ik}(t)  \int_{t-\tau_{ik}(t)}^t\ga_{ik}(s)\hat h_{ik}(s,\hat x_i(s))\,  d_s\eta_{ik}(t,s),\  i=1,\dots,n,
\end{split}
\end{equation*}
where $ \hat a_{ij}(t)= v_i^{-1}a_{ij}(t) v_j, j\ne i$, and $\hat h_{ik}(t,x)=v_i^{-1}h_{ik}(t,v_ix)$. The matrix $D(t)-\big [\hat a_{ij}(t)\big]$ still satisfies (H2). In this way, and dropping the hats for simplicity, we consider the original system \eqref{0}, but   suppose that (H5) holds with $v={\bf 1}:=(1,\dots,1)$.
Thus, there exist constants $\eta_i>0\, ( i=1,\dots,n)$ such that
$$\be_i(t)\ge d_i(t)-\sum_{j\ne i}a_{ij}(t)+\eta_i,\q t\in\R.$$
On the other hand, $d_i(t)-\sum_{j\ne i}a_{ij}(t)\le \ol d_i:=\max _{t\in [0,\omega]}d_i(t)$, and with $1<\al_i<1+\eta_i/\ol d_i$ we obtain
\begin{equation}\label{2.7}
\al_i^{-1}\be_i(t)-d_i(t)+\sum_{j\ne i}a_{ij}(t)>0,\q {\rm for}\q t\in\R,  i=1,\dots,n.
\end{equation}

  From the dissipativeness of the system asserted in (ii), and for $h_i^-$ as in (H4),
 we can choose
  $L>m>0$  such that the uniform estimate
  \begin{equation}\label{2.8}
\limsup_{t\to\infty} x_i(t,0,\phi)< L\q {\rm for}\q \phi\in C_0,\ i=1,\dots,n,
\end{equation} 
holds and $h_i^-(m)=\min_{x\in [m,L]}h_i^-(x)$, with
  $(h_i^-)'(x)>0$ and $\al_i^{-1}x<h_i^-(x)$ for $x\in(0,m]$ and  all $i$.

Consider the auxiliary  {\it cooperative} system
  \begin{equation}\label{0H}
  \begin{split}
x_i'(t)=&-d_i(t)x_i(t)+\sum_{j=1,j\ne i}^n a_{ij}(t)x_j(t)\\
&+\sum_{k=1}^m \be_{ik}(t)  \int_{t-\tau_{ik}(t)}^t\ga_{ik}(s)H_i(x_i(s))\,  d_s\eta_{ik}(t,s)=:F_i(t,x_t),\  i=1,\dots,n,
\end{split}
\end{equation}
where $H_i(x)=h_i^-(x)$ if $0\le x\le m$, $H_i(x)=h_i^-(m)$ if $x\ge m$. 

For $x(t)$ a positive solution of \eqref{0},   for $t>0$ sufficiently large and $1\le i\le n$, we have $x_i(t)\le L$ and $h_{ik}(t,x_i(t))\ge H_i(x_i(t))$. Therefore, if \eqref {0H} is uniformly persistent, then \eqref{0} is uniformly persistent as well \cite{Smith}. 

Now,  we consider any solution $x(t)=x(t,\sigma,\phi,F)$  of \eqref {0H} with  $x_{\sigma}=\phi\in C_0$ and $\sigma\in\R$. We claim that there is $T=T(\sigma,\phi)\ge \sigma$ such that
 \begin{equation}\label{2.9}
  x_i(t)\ge m\q {\rm for}\q  t\ge T,1\le i\le n.
\end{equation}

\smal

We first  prove that if $\displaystyle \min\{ x_j(t):1\le j\le n,t\in [T,T+\tau]\}\ge m $ for some $T\ge \sigma$, then $ x_j(t)\ge m $ for all $t\ge T$ and $j=1,\dots, n$. 

For simplicity of exposition, take $T=\sigma=0$.
Assume that $x_j(t)\ge m$ for $t\in [0,\tau]$ and $j=1,\dots, n$. Let $t_0\in [\tau, 2\tau]$ and $i\in \{1,\dots, n\}$ be such that $x_i(t_0)=\min\{ x_j(t):1\le j\le n,t\in [\tau, 2\tau]\}$.
We have 
$$0\ge x_i'(t_0)=-d_i(t_0)x_i(t_0)+\sum_{j\ne i} a_{ij} (t_0)x_j(t_0)+\sum_{k=1}^m \be_{ik} (t_0) \int_{t_0-\tau_{ik}(t_0)}^{t_0} \!\! \ga_{ik}(s)H_i(x_i(s))\,  d_s\eta_{ik}(t_0,s). 
 $$
Suppose that  $x_i(t_0)<m$. For $k=1,\dots,m$ and  $s\in[t_0-\tau_{ik}(t_0),t_0]\subset [0,t_0]$, we have $x_i(s)\ge x_i(t_0)$, hence $H_i(x_i(s)) \ge H_i(x_i(t_0))$. From  the definition of $\be_i(t)$ in \eqref{beta-i}, we obtain 
  \begin{equation}\label{06}
\begin{split}
0&\ge \left (-d_i(t_0)+\sum_{j=1}^n a_{ij}(t_0)\right ) x_i(t_0) +\be_i(t_0)H_i(x_i(t_0)) \\
&\ge \left (-d_i(t_0)+\sum_{j=1}^n a_{ij}(t_0)+\al_i^{-1}\be_i(t_0)\right )x_i(t_0)>0,
\end{split}
\end{equation}
which is not possible. Thus, $x_i(t_0)\ge m$. This implies that $x_j(t)\ge m$ on $[0,2\tau]$ for all $j=1,\dots,n$. By iteration, we obtain the same lower bound $m$ on $[0,\infty)$, which proves  \eqref{2.9}.

Next, we need to show that there exists an interval of length $\tau$ where the minima of all components $x_i(t)$ are larger to or equal to $m$. The proof  follows  by  adapting slightly the arguments in \cite{FOS}, so we do not include it  here.\qed\end{pf}

\begin{rmk} For systems \eqref{3} with all the coefficients and delay functions continuous and bounded, if we replace (H2) and (H5) by slightly stronger assumptions,  the claims in the above theorem remain valid without assuming that the coefficients and delay functions are periodic. See \cite{FOS} for details, as well as for supplementary results. See also \cite{ObayaSanz}, for the uniform persistence of a Nicholson almost-periodic system with one constant delay in each equation.
\end{rmk}

\section{Existence of a positive periodic solution}
\setcounter{equation}{0}


In the case of periodic systems,   we now show that the criterion for uniform persistence in Theorem \ref{thm2.1}(iii) also provides a criterion for the existence of a positive $\omega$-periodic solution.  We start with some algebraic definitions.

\begin{defn}\label{def3.1} Let $A=[a_{ij}]$ be a square matrix. We say that $A$ is {\it nonnegative}, and write $A\ge 0$, if all its entries are nonnegative. For $A$  with nonpositive off-diagonal entries   (i.e., $a_{ij}\le 0$ for $i\ne j$),
 $A$ is said to be a {\it non-singular M-matrix} or a {\it matrix of class K}
if   all its eigenvalues have positive  real parts. 
\end{defn}

We remark that some authors use simply the term {\it M-matrix} to designate a  {\it non-singular M-matrix}.
There are many alternative equivalent definitions of  non-singular  M-matrices, see e.g.~\cite{Fiedler}. 
Namely,   for a square matrix $A$ with nonpositive off-diagonal entries, the following conditions are equivalent: 
(i) $A$  is a  non-singular M-matrix;
(ii) there exists a vector $u\gg 0$ such that $Au\gg 0$; 
(iii) $A$ is non-singular and $A^{-1}\ge 0$.

 \begin{thm}\label{thm3.1}   Assume  (H0)-(H5).  Then \eqref{0} has a positive $\omega$-periodic  solution.
\end{thm}

\begin{pf} The proof will be divided in several steps.

\smal

(i) From Theorem \ref{thm2.1}(i), the linear homogeneous ODE \eqref{2.1}
is exponentially asymptotically stable. Let $K\ge 1,\al>0$ be such that $|X(t)X^{-1}(s)|\le Ke^{-\al(t-s)}$ for $t\ge s$, where $X(t)$ is the fundamental matrix solution of \eqref{2.1} with $X(0)=I$.
For $y_0\in\R^n$, the solution of \eqref{2.1} with initial condition $y(s)=y_0$ is given by $X(t)X^{-1}(s)y_0$. It was observed in Section 2 that \eqref{2.1} is cooperative, hence  its solutions are monotone relative to the order  in $\R^n$, i.e., $y(t,s,x_0)\le y(t,s,y_0)$  if $x_0\le y_0$.
Moreover, for $s\in\R,y_0\in\R^n$, a solution $y(t)=y(t,s,y_0)$ of  \eqref{2.1} satisfies
$y_i'(t)\ge -d_i(t)y_i(t),\, i=1,\dots,n$, and therefore $y(t,s,y_0)\ge  0$ for $t\ge s$ whenever $y_0\ge 0$, with $y_i(t,s,y_0)=(X(t)X^{-1}(s)y_0)_i>0$ if  $y_{0i}>0$, for any $1\le i\le n$. For $t\ge s$, we derive that $X(t)X^{-1}(s)\ge 0$, i.e., all entries of the matrices $X(t)X^{-1}(s)$ are nonnegative  (and that their   diagonal entries are positive). For the monodromy matrix $C=X(\omega)$, we have $C=X^{-1}(t)X(\omega+t)$ for $t\in\R$.
 The matrices $T(t):=X(\omega+t)X^{-1}(t),\,  t\in\R,$  are nonnegative, $\omega$-periodic and  similar to $C$. Since all the characteristic multipliers of \eqref{2.1}  have moduli less than one, the spectral radius $\rho (T(t))$ of  $T(t)$ is less than one.
 Consequently, $I-T(t)$ is a nonsingular M-matrix, and has  inverse $(I-T(t))^{-1}\ge 0$.

\med

(ii) From (H5), there exists  $v=(v_1,\dots,v_n)\gg 0$ such that
\begin{equation}\label{5.1}
\eta_i:=\min_{t\in [0,\omega]}\Big(\be_i(t)v_i-d_i(t)v_i+\sum_{j\ne i}a_{ij}(t)v_j\Big)>0, \q i=1,\dots,n.
\end{equation}

As before, we effect the scaling of variables    $\hat x_j(t)=x_j(t)/v_j$ in \eqref{0}, and obtain a new system of the form \eqref{0}  where $ \hat a_{ij}(t)= v_i^{-1}a_{ij}(t) v_j, j\ne i$, and $\hat h_{ik}(t,x)=v_i^{-1}h_{ik}(t,v_ix)$.
 Hence, and without loss of generality, we may consider   \eqref{0}
and  take   $v=(1,\dots,1)={\bf 1}$ in (H5).  
 As in the proof of Theorem \ref{thm2.1}, we deduce that there are
 constants $\al_i>1$ such that \eqref{2.7} is satisfied.
 
  Theorem \ref{thm2.1}(iii) implies that \eqref{0} is permanent.  
 Consider    uniform lower and upper bounds $m,L$ for all positive solutions of \eqref {0}, as in the uniform estimates \eqref{3.5}.
As a result of (H4), we can choose 
  $L>m>0$, with $m$ sufficiently small such that $h_i^-(m)=\min_{x\in [m,L]}h_i^-(x)$,  with $h_i^-(x)$ increasing on $[0,m]$ and
   $$\al_i^{-1}x<h_i^-(x)\q {\rm for}\q x\in(0,m],\q i=1,\dots,n.$$

\med

(iii) We  introduce some further terminology. For simplicity, we may suppose that $\tau\ge \omega$  
(otherwise we choose some $\bar \tau\ge \omega$ and insert $C$ into $C([-\bar\tau,0];\R^n$)). For $\phi\in C$ such that $\phi(t+\omega)=\phi(t)$ for all $t,t+\omega\in [-\tau,0]$, we write $\tilde \phi$ for the $\omega$-periodic function  defined in $\R$ which coincides with $\phi$ on $[-\tau,0]$.
 Denote by $C_\omega, C_\omega^+$  the sets of $\omega$-periodic continuous functions $\phi:\R\to \R^n$, respectively $\phi:\R\to \R_+^n$, which can be identified as  subsets of $C, C^+$, respectively, with the same topology. 
 
 Now, suppose that $x(t)=x(t,-\tau,\phi)$ is a solution of \eqref{0}, with initial condition $x_{-\tau}=\phi\in C_\omega$.
 By the  variation of constants formula for ODEs, 
\begin{equation}\label{VCF}
x(t)=X(t)X^{-1}(t_0)x(t_0)+X(t)\Big (\int_{t_0}^t X^{-1}(s) M(s,x_s)\, ds\Big)\q (t,t_0\ge -\tau),
\end{equation}
where $M(t,\phi)=(M_1(t,\phi_1),\dots, M_n(t,\phi_n))$ is given by \eqref{M0} and, as before, $x_s=x_{|_{[s-\tau,s]}}$ for $s\ge -\tau$.
Clearly, $x(t)$ is $\omega$-periodic if and only if   
$x_{\omega}=x_0$.
From \eqref{VCF}, $x(\omega+\th)=x(\th)$ for  $ \th\in [-\tau,0]$  if and only if 
\begin{equation}\label{5.3}
x(\th)=X(\omega+\th)X^{-1}(\th)x(\th)
+X(\omega+\th)\int_{\th}^{\omega+\th}X^{-1}(s)M (s, x_s )\, ds,
 \end{equation}
for $ \th\in [-\tau,0]$.
 This is equivalent to saying that 
$x(t)$   is a fixed point of  the operator ${\cal F}: C_\omega \to C$ defined by
\begin{equation}\label{calF}
 ({\cal F}\phi)(\th)=\Big (I- T(\th)\Big )^{-1} \left (X(\omega+\th)\int_\th^{\omega+\th}X^{-1}(s)M(s, \phi_s)\, ds\right ),\ \phi\in C_\omega,\th\in  [-\tau,0].
\end{equation}
Thus, we look for a fixed point $\phi\in\, int(C_\omega^+)$  of  the operator ${\cal F}$.

(iv) The aim is to apply the Schauder fixed point theorem to the operator ${\cal F}$ in an appropriate  subset of $int(C_\omega^+)$.

We first claim that ${\cal F}\phi\in C_\omega^+$ whenever $\phi\in C_\omega^+$. 
Let $\phi\in C_\omega$, and set
\begin{equation}\label {G}
G(t;\phi):=X(\omega+t)\int_t^{\omega +t}X^{-1}(s)M(s, \phi_s)\, ds,\q t\in\R.
\end{equation} 
We have
$$
G(\omega+t;\phi)=\int_t^{\omega +t} X(2\omega+t) X^{-1}(\omega+s)M(\omega+s,\phi_{\omega+s})\, ds
=G(t;\phi),\q t\in\R,$$
because $ \phi(t)$ and $t\mapsto M(t,\psi)$ are  $\omega$-periodic  and $X(2\omega+t) X^{-1}(\omega+s)=X(\omega +t)CX^{-1}(\omega+s)=X(\omega +t)X^{-1}(s).$ From step (i), $T(t)$ is also $\omega$-periodic, hence ${\cal F}\phi\in C_\omega$. Since  $X(\omega +t)X^{-1}(s)\ge 0, (I-T(\th))^{-1}\ge 0$, we further  derive that  ${\cal F}\phi\ge 0$ for $\phi\in C_\omega^+$.
  
 From the continuity of $(I-T(\th))^{-1}$,
there exists $c=\max_{\th\in[-\th,0]} |(I-T(\th))^{-1}|$, thus 
 $$\|  {\cal F}\phi\|\le cK\be^0L^0 \frac{1}{\al} (1-e^{-\al \omega}),$$
 where  $L^0,\be^0$ are such that $h_{ik}(t,x)\le L^0$ and  $\be_i(t)\le \be^0$ for all $i,k$ and $t\in\R,x\ge 0$.  Therefore,  ${\cal F}$ transforms  $C_\omega^+$ into a bounded set of $C_\omega^+$. Choose $R\ge L$, for $L$ as in \eqref{3.5}, such that ${\cal F}(C_\omega^+)\subset [0,R{\bf 1}]_\omega$, where $[0,R{\bf 1}]_\omega:=\{\phi\in C_\omega^+:\phi_i\le R, 1\le i\le n\}$.

 We now prove that ${\cal F}(C_\omega^+)$ is equicontinuous in $C_\omega^+$. 
For $\phi\in C_\omega^+$, consider $G(t;\phi)$ as in \eqref{G}, and
 $t_1,t_2\in [-\omega,0]$. We have
\begin{equation}\label{3.6}
 \begin{split}
 |({\cal F}\phi)(t_1)-({\cal F}\phi)(t_2)|&\le c\, |G(t_1;\phi)-G(t_2;\phi)|\\
 &+\left |  \big (I- T(t_1)\big )^{-1}-\big (I- T(t_2)\big )^{-1}\right |G(t_2;\phi)\, .\\
\end{split}
 \end{equation}
Observe that
 \begin{equation*}
 \begin{split}
& \left | \int_t^{\omega +t}X^{-1}(s) M(s,\phi_s)\, ds\right |\le  K\be^0L^0 \frac{1}{\al} (e^{\al \omega}-1)=:C_1,\\
&|G(t;\phi)|\le K\be^0L^0 \frac{1}{\al} (1-e^{-\al \omega}),\q \forall \phi\in C_\omega^+,\forall t\in [-\omega,0],
\end{split}
 \end{equation*}
 and 
\begin{equation*}
 \begin{split}
 |G(t_1;\phi)-G(t_2;\phi)|&\le \left |X(t_1+\omega)-X(t_2+\omega)\right |C_1\\
 &+K\left | \int_{t_1}^{t_2}X^{-1}(s) M(s,\phi_s)\, ds-\int_{\omega+t_1}^{\omega+t_2}X^{-1}(s) M(s,\phi_s)\, ds \right |\\
 &\le \left |X(t_1+\omega)-X(t_2+\omega)\right |C_1+ K\be^0L^0 \frac{1}{\al} (1+e^{-\al \omega})\Big |e^{-\al t_1}-e^{-\al t_2}\big |, 
 \end{split}
 \end{equation*}
for all $\phi\in C_\omega^+,t_1,t_2\in [-\omega,0].$
Inserting these estimates in \eqref{3.6}, we conclude that the family ${\cal F}(C_\omega^+)$ is equicontinuous. By Ascoli-Arzel\`a theorem, ${\cal F}(C_\omega^+)$  is relative compact in $C_\omega^+$.

Next, we   claim that
 \begin{equation}\label{5.6}{\cal F}([m{\bf 1},\infty)_\omega)\subset [m{\bf 1},\infty)_\omega,
 \end{equation}
 where $[m{\bf 1},\infty)_\omega:=\{\phi\in C_\omega^+:\phi_i\ge m, 1\le i\le n\}$. Note that all solutions $x(t)=x(t,\sigma,\phi)$ of \eqref{0} (with $\phi\in C_0$) satisfy $m\le x_i(t)\le R$ for $t$ sufficiently large and $1\le i\le n$.
 
Take  $\phi\in C_\omega$ with $\phi(s)\ge m{\bf 1}$ for $s\in \R$. From step (ii),
we have  $h_{ik}(s,\phi_i(s))\ge h_i^-(\phi_i(s))\ge \al_i^{-1} \phi_i(s)\ge \al_i^{-1}m$ for all $i,k$ and $s\in\R$. For $M_i$ defined in \eqref{M0}, we obtain $M_i(s,\phi_{i,s})\ge   \be_i(s)\al_i^{-1} m$, and 
 \eqref {2.7} yields $M_i(s,\phi_{i,s})\ge  m [d_i(s)-\sum_{j\ne i} a_{ij}(s)], i=1,\dots,n, s\in\R$. Since $X(\omega +\th)X^{-1}(s)\ge 0$, we deduce that
  \begin{equation}\label{5.7}
  X(\omega+\th)\int_\th^{\omega+\th}X^{-1}(s)M(s,\phi_s)\, ds\ge m X(\omega+\th)\int_\th^{\omega+\th}X^{-1}(s)[D(s)-A(s)]{\bf 1}\, ds.
  \end{equation}
  The differentiation of the identity $I=X^{-1}(s)X(s)$ leads to 
  $\frac{d}{ds} \Big (X^{-1}(s)\Big )=X^{-1}(s)[D(s)-A(s)].$
 From \eqref{5.7},  we derive
  $$X(\omega+\th)\int_\th^{\omega+\th}X^{-1}(s)M(s,\phi_s)\, ds\ge  m [I- T(\th)\Big ]{\bf 1},$$
  and finally from \eqref{calF}  obtain
  $$
 ({\cal F}\phi)(\th)\ge m \Big [I- T(\th)\Big ]^{-1} \Big [I- T(\th)\Big ]{\bf 1}=m{\bf 1},$$
 which proves the claim \eqref{5.6}.

Consider the  convex, closed bounded subset $[m{\bf 1} ,R{\bf 1}]_\omega:=\{\phi\in C_\omega: m{\bf 1}\le \phi \le R{\bf 1} \}$ of $C_\omega$.  Applying  Schauder's fixed point theorem to the restriction  (still denoted by  ${\cal F}$) ${\cal F}:[m{\bf 1} ,R{\bf 1}]_\omega\to [m {\bf 1} ,R{\bf 1}]_\omega$, we conclude that there exists a fixed point $\phi^*\in [m {\bf 1} ,R{\bf 1}]_\omega$.  From (ii),  $\phi^*(t)$ is an $\omega$-periodic solution of \eqref{0}. The proof is complete.
\qed\end{pf}

A general criterion concerning the existence of a positive periodic solution for periodic $n$-dimensional Nicholson systems  is trivially  obtained as a consequence of Theorem \ref{thm3.1}.
 
\begin{thm}\label{thm3.2}   Consider  \eqref{1.2'}  where all the functions $d_i(t),a_{ij}(t), \be_{ik}(t), \ga_{ik}(t),c_{ik}(t), \tau_{ik}(t)$ satisfy (H0)-(H3) and (H5). 
Then there exists (at least) one positive $\omega$-periodic solution of \eqref{1.2'}. A similar results holds for \eqref{1.2}, with $\be_i(t)$ in \eqref{beta-i} replaced by $\be_i(t)=\sum_{k=1}^m \be_{ik}(t)$.\end{thm}

As a by-product, Theorem \ref{thm3.1} also provides conditions for  the existence of a positive equilibrium  for systems with {\it autonomous} coefficients.

\begin{thm}\label{thm3.3} Consider the system
 \begin{equation}\label{3.14}
x_i'(t)=-d_ix_i(t)+\sum_{j=1,j\ne i}^n a_{ij}x_j(t)+ \sum_{k=1}^m \be_{ik} h_{ik}(x_i(t-\tau_{ik}(t))), \  i=1,\dots,n,t\ge 0,
\end{equation}
where $d_i>0$, $a_{ij}\ge 0, \be_{ik}\ge 0$ with $\be_i:=\sum_{k=1}^m \be_{ik}>0$,
$ \tau_{ik}:[0,\infty)\to [0,\infty)$  are continuous and uniformly bounded from above by some $\tau>0$, and  
\begin{itemize}
\item[(H4*)]   $h_{ik}:[0,\infty)\to [0,\infty)$ are bounded,  locally Lipschitz continuous on $[0,\infty)$ and continuously differentiable on a right  neighborhood of $0$,  with $h_{ik}(0)=0,h_{ik}'(0)=1$ and $ h_{ik}(x)>0$ for $x>0$,   
\end{itemize}
for all $i,j=1,\dots,n, k=1,\dots,m$. Define the $n\times n$ matrices
\begin{equation}\label{3.15}
A=[a_{ij}],\ B=diag\, (\be_1,\dots,\be_n),\ D=diag\, (d_1,\dots,d_n),\ M=B-D+A,
\end{equation}
where $a_{ii}:=0\ (1\le i\le n)$.  Assume that:
(i) $D-A$ is a non-singular M-matrix;
(ii) $Mv\gg 0$ for some vector $v\gg 0$. 
  Then \eqref{3.14} has a positive equilibrium.
\end{thm}

\begin{pf}
For $D,A$  as in \eqref{3.15}, hypotheses (H1),~(H2) are satisfied, thus the linear autonomous ODE $x'=-(D-A)x$ is exponentially asymptotically stable. Together with \eqref{3.14}, consider its associated ODE  without delays:
\begin{equation}\label{5.8'}
x_i'(t)=-d_ix_i(t)+\sum_{j=1,j\ne i}^n a_{ij}x_j(t)+\sum_{k=1}^m \be_{ik}    h_{ik}(x_i(t)),\
i=1,\dots,n,\ t\ge 0.
\end{equation}
Systems \eqref{3.14} and \eqref{5.8'} have the same equilibria. We apply Theorem \ref{thm3.1}  to \eqref{5.8'}, noticing that in this case  $\omega=0,\tau=0$ and $C=\R^n$, and deduce the existence of a positive equilibrium.
\qed\end{pf}

\begin{rmk}\label{rmk3.1}  Consider the case of  autonomous ODEs $x'=f(x)$ with $f:\R^n\to\R^n$  a $C^1$ function. From the work of Hofbauer  \cite{Hofbauer}, it follows that if $x'=f(x)$ is dissipative and the nonnegative cone $[0,\infty)^n$ is forward invariant for its flow, then there exists a {\it saturated equilibrium} $x^*\ge 0$ (see \cite{Hofbauer}  for a definition), which may however lay on the border of $[0,\infty)^n$. Supplementary results can be found in  \cite{Hale,Hofbauer}. On the other hand, if in addition $f(0)=0$ and  the system $x'=f(x)$ is uniformly persistent in $[0,\infty)^n\setminus \{ 0\}$, obviously there are no nonnegative  equilibria $x^*$ besides the trivial one, hence  a positive equilibrium must exist. For the ODE \eqref{5.8'},
Theorem \ref{thm3.3} asserts the existence of such an equilibrium without demanding the $C^1$-smootheness of the  vector field, though; in fact, from our assumptions the vector field in \eqref{5.8'} is simply  locally Lipschitzian (in order to guarantee the uniqueness of solutions) and continuously differentiable in a vicinity of $0^+$. Moreover, for the case of  {\it autonomous} Nicholson systems \eqref{1.2} (thus  with constant coefficients  and delays) with  $c_{ik}(t)\equiv 1$,
 the existence of a positive equilibrium was established in \cite{FariaRost} exactly under the conditions in Theorem \ref{thm3.3}. 

\end{rmk}

From  Theorem \ref{thm3.2}, we recover or improve some results in the  literature. 

\begin{cor}\label{cor3.1} Consider the equation
\begin{equation}\label{n=1}
 x'(t)=-d(t)x(t)+\sum_{k=1}^m\be_k(t) h_k(t,x(t-\tau_k(t))),
 \end{equation} 
where the functions $d(t),\be_{k}(t),\tau_k(t)$ are continuous, non-negative and $\omega$-periodic, with $d(t)>0$ for $t\in\R$, and $h_k(t,x)$ satisfy (H4). If 
\begin{equation}\label{1.5}
\sum_{k=1}^m \be_k(t)>d(t),\q t\in [0,\omega],
 \end{equation}  
  then  there exists a positive $\omega$-periodic solution of \eqref{n=1}. 
 \end{cor}

 
 \begin{cor}\label{cor3.2} Consider the periodic Nicholson's equations with distribute delays
\begin{equation}\label{NPP-dist}
 x'(t)=-d(t)x(t)+\sum_{k=1}^m\be_k(t)  \int_{t-\tau_{k}(t)}^t\ga_{k}(s)x(s)e^{-c_k(s)x(s)}\,  ds,
 \end{equation} 
where  $d(t), c_{k}(t)>0,\be_{k}(t),\ga_k(t),\tau_k(t)\ge 0$ are continuous and  $\omega$-periodic. If
$$
\sum_{k=1}^m \Big (\be_k(t)\int_{t-\tau_{k}(t)}^t\ga_{k}(s)\, ds\Big )>d(t),\q t\in [0,\omega],
$$
 then \eqref{NPP-dist}   has a positive $\omega$-periodic solution.  In particular, for
 the equation
\begin{equation}\label{NPP-dist m=1}
x'(t)=-d(t)x(t)+\be(t)  \int_{t-\tau(t)}^t\ga(s)x(s)e^{-c(s)x(s)}\,  ds,
 \end{equation}
 there is a positive $\omega$-periodic solution  if $\be(t)\int_{t-\tau(t)}^t\ga(s)\, ds>d(t),\, t\in [0,\omega]$. 
 \end{cor}

 \begin{rmk}\label{rmk3.1'}   For the periodic Nicholson's equations with multiple discrete delays,
\begin{equation}\label{NPP}
 N'(t)=-d(t)N(t)+\sum_{k=1}^m\be_k(t) N(t-\tau_k(t))e^{-c_k(t)N(t-\tau_k(t))},
  \end{equation} 
where  $d(t)>0, c_{k}(t)>0,\be_{k}(t)\ge 0,\tau_k(t)\ge 0$ are continuous and  $\omega$-periodic,  Li and Du \cite{LiDu} proved the existence of a positive $\omega$-periodic solution if \eqref{1.5} holds,  by using the Krasnoselskii  fixed point theorem on cones. It is clear that the result in \cite{LiDu} follows as a particular case of
 Corollary \ref{cor3.1}. 
 On the other hand, Corollary \ref{cor3.2}  improves the result in \cite{AmsterIdels}, where the existence of a  positive $\omega$-periodic solution for \eqref{NPP-dist m=1} was obtained  under the stronger condition $$\min_{t\in [0,\omega]}\ga(t)>\max_{t\in [0,\omega]}\frac{d(t)}{\tau(t)\be(t)}.$$
 \end{rmk}

 For $n=2$, the hypotheses (H2), (H5) are also easily verifiable in practice. 
 For illustration, we state here a criterion for  systems with discrete time-varying delays.
 
 \begin{cor}\label{cor3.3} Consider the planar  system given by
 \begin{equation}\label{3.160}
 \begin{split}
x_1'(t)&=-d_1(t)x_1(t)+ a_1(t)x_2(t)+\sum_{k=1}^{m_1} \be_{1k}(t)    h_{1k}(t,x_1(t-\tau_{1k}(t)))\\
x_2'(t)&=-d_2(t)x_2(t)+ a_2(t)x_1(t)+\sum_{k=1}^{m_2} \be_{2k}(t)     h_{2k}(t,x_2(t-\tau_{2k}(t)))
\end{split}
\end{equation}
where  $m_1,m_2\in \N$, $d_i(t), a_i(t), \be_{ik}(t),t\mapsto h_{ik}(t,x)\, (x\ge 0),\tau_{ik}(t)$ are continuous, nonnegative  and $\omega$-periodic, with $d_i(t), a_i(t)$ and $\be_i(t):=\sum_{k=1}^{m_i} \be_{ik}(t)$ strictly positive for $t\in [0,\omega]$, and $h_{ik}$ satisfy (H4), 
 $i=1,2,k=1,\dots ,m_i$. In addition, suppose that:\vskip 0cm
(i) $\dps \min_{t\in [0,\omega]} \frac{d_1(t)}{a_1(t)}> \max_{t\in [0,\omega]} \frac{a_2(t)}{d_2(t)}$;
\vskip .1cm
(ii) there exist constants $u_1,u_2>0$ such that
$$u_1(\be_1(t)-d_1(t))+u_2a_1(t)>0,\ u_2(\be_2(t)-d_2(t))+u_1a_2(t)>0,\q t\in [0,\omega].$$
Then  \eqref{3.160} has at least one  positive $\omega$-periodic solution.\end{cor}

\begin{pf}  From condition (i), choose  $v_2$ with 
$\max_{t\in [0,\omega]} \frac{a_2(t)}{d_2(t)}<v_2<\min_{t\in [0,\omega]}\frac{d_1(t)}{a_1(t)}.$
With $v=(1,v_2)$, we have
$$\left[ \begin{array}{cc} d_1(t)&-a_1(t) \\
-a_2(t)&d_2(t)
\end{array}\right]v>0,\q t\in [0,\omega],$$
thus (H2) is satisfied. On the other hand, (ii) is hypothesis (H5) for the case $n=2$, and the result follows from Theorem \ref{thm3.1}.
 \qed \end{pf}

 \begin{rmk}\label{rmk3.2}  Liu \cite{Liu11} considered the planar Nicholson system 
  \begin{equation}\label{3.16}
 \begin{split}
x_1'(t)&=-d_1(t)x_1(t)+ a_1(t)x_2(t)+\sum_{k=1}^{m_1} \be_{1k}(t)    x_1(t-\tau_{1k}(t))e^{-c_{1k}(t)x_1(t-\tau_{1k}(t))}\\
x_2'(t)&=-d_2(t)x_2(t)+ a_2(t)x_1(t)+\sum_{k=1}^{m_2} \be_{2k}(t)    x_2(t-\tau_{2k}(t))e^{-c_{2k}(t)x_1(t-\tau_{2k}(t))}
\end{split}
\end{equation} 
with all coefficients and delay functions  $\omega$-periodic, continuous and positive. By constructing a suitable Lyapunov functional, Liu obtained the existence (and uniqueness) of a positive $\omega$-periodic solution  by imposing some   other rather restrictive constraints. Among these additional conditions, it was assumed that (cf.~\cite[Theorem 2.1]{Liu11})
$$ \min_{t\in [0,\omega]} \left ( \sum_{k=1}^{m_i} \be_{ik}(t)-d_i(t)\right)>0,\
  \max_{t\in [0,\omega]} a_i(t) +e^{-2}\sum_{k=1}^{m_i}  \max_{t\in [0,\omega]}\be_{ik}(t)<\min_{t\in [0,\omega]}  d_i(t),\, i=1,2,$$ which are assumptions stronger than (i),(ii) in Corollary \ref{cor3.3}.
 \end{rmk}

\begin{rmk}\label{rmk3.3} As observed in the Introduction, our results are not optimal, and better criteria involving the average of the periodic coefficients $d_i(t),a_{ij}(t),\be_i(t)$ in \eqref{0} are desirable. In fact, even for the case of $n=1$ with one discrete delay, our method does not allow to recover the criterion of Chen \cite{Chen},  who establish the existence  of a positive $\omega$-periodic solution of \eqref{NP}
  under the conditions
\begin{equation*}
\begin{split} 
\ol \be> \ol d \exp(2\omega \ol d),&\q {\rm if}\q \tau(t)\ {\rm is}\ \omega{\rm -periodic}\\
\ol \be> \ol d,&\q {\rm if}\q \tau(t)=m\omega
  \end{split} 
   \end{equation*} 
 where $\ol \be:=\frac{1}{\omega}\int_0^\omega \be(t)\,dt,  \ol d:= \frac{1}{\omega}\int_0^\omega d(t)\, dt,$
and $m$ is some positive intege. Another limitation of our approach is that it cannot be applied directly when there exist some $i\in\{1,\dots,n\}, t_0\in [0,\omega]$ such that either $d_i(t_0)=0$ or $\be_i(t_0)$ (see \cite{FLT} for an example).\end{rmk}

\begin{exmp}\label{exmp3.1} Consider the $\pi$-periodic planar system of Mackey-Glass type 
 \begin{equation}\label{3.17}
 \begin{split}
x_1'(t)&=-(\epsilon_1+\sin^2t)x_1(t)+ |\cos (2t)|x_2(t)+   \frac{(\de_1+\cos^2t)x_1(t-\sin^2t)}{1+e^{-\sin^2t} x_1^\al(t-\sin^2t)}\\
x_2'(t)&=-(\epsilon_2+\cos^2t)x_2(t)+|\cos (2t)|x_1(t)+    \frac{(\de_2+\sin^2t)x_2(t-\cos^2t)}{1+(2+\cos (2t))x_2^\be(t-\cos^2t)}
\end{split}
\end{equation}
where  $\epsilon_i,\de_i>0$ for $i=1,2$ and $\al,\be\ge 1$. The nonlinearities have the form \eqref{hiMG}. For  \eqref{3.17} and with the notation in \eqref{2.3},
 \begin{equation*}
\begin{split}D(t)-A(t)&=\left[ \begin{array}{cc} \epsilon_1+\sin^2t&-|\cos (2t)| \\
-|\cos (2t)|&\epsilon_2+\cos^2t
\end{array}\right],\\
M(t)&=\left[ \begin{array}{cc} (\de_1+\cos^2t)-(\epsilon_1+\sin^2t)&|\cos (2t)| \\
|\cos (2t)|&(\de_2+\sin^2t)-(\epsilon_2+\cos^2t)
\end{array}\right].
\end{split}
\end{equation*}

In addition, suppose that 
$\epsilon_1\epsilon_2\ge  1, \, \de_1>\epsilon_1, \, \de_2>\epsilon_2.$ 
Note that
 $$M(t)\left[ \begin{array}{c}1\\1\end{array}\right]=\left[ \begin{array}{c} \de_1-\epsilon_1+\cos(2t)+|\cos(2t)|\
\\ \de_2-\epsilon_2-\cos(2t)+|\cos(2t)|\end{array}\right].
$$ 
Clearly, (H5) is satisfied with $u=(1,1)$. 
Next, take $v_2>0$ such that
$\epsilon_2^{-1}\le v_2\le \epsilon_1.$
With $v=(1,v_2)$, we obtain 
$$[D(t)-A(t)]v\ge  \left[ \begin{array}{c} \sin^2 t\\ \epsilon_2^{-1}\cos^2t\end{array}\right]\q  {\rm for}\q t\in [0,\pi].$$
Thus, assumption (H2) holds. From Theorem \ref{thm3.1}, we conclude that \eqref{3.17} has a $\pi$-periodic positive solution.
\end{exmp}

\begin{exmp}\label{exmp3.2} Consider the planar Nicholson system
 \begin{equation}\label{3.18}
  \begin{split}
x_1'(t)&=-(\epsilon_1+\cos^2 t)x_1(t)+ a_{12}e^{-2+\sin^2t}x_2(t)+e^{\cos^2 t}
\int_{t- (\be_1e^{-\cos^2 t}+1)}^t \!\! x_1(s)e^{-(1+|\sin s|)x_1(s)}\, ds\\
x_2'(t)&=-(\epsilon_2+\sin^2t)x_2(t)+a_{21}e^{\cos^2t}x_1(t)+e^{\sin^2 t}
 \int_{t-(\be_2e^{-\sin^2 t}+1)}^t \!\! x_2(s) e^{-e^{\sin(2s)}x_2(s)}\, ds,
\end{split}
\end{equation}
where  $a_{12}, a_{21},\epsilon_i,\be_i>0$ for $i=1,2$. The functions $\be_i(t)$  in \eqref{beta-i} are given by
$$\be_1(t)=e^{\cos^2 t}
\int_{t- (\be_1e^{-\cos^2 t}+1)}^t \!\!ds=\be_1+e^{\cos^2 t},\q \be_2(t)=e^{\sin^2 t}
 \int_{t-(\be_2e^{-\sin^2 t}+1)}^t \!\! ds=\be_2+e^{\sin^2 t},$$ for $t\in\R$.  For $y=\cos^2t$ and   matrices defined as in \eqref{2.3}, we have
 \begin{equation*}
D(t)-A(t)=\left[ \begin{array}{cc} \epsilon_1+y&-a_{12}e^{-(y+1)} \\
-a_{21}e^y&\epsilon_2+1-y
\end{array}\right],M(t)=\left[ \begin{array}{cc} \be_1+e^y&0 \\
0& \be_2+e^{1-y}
\end{array}\right]+A(t)-D(t).
\end{equation*}

Suppose that
 \begin{equation}\label{3.19}
 \epsilon_1\epsilon_2+\min\{\epsilon_1,\epsilon_2\}\ge a_{12}a_{21}.
 \end{equation} In this case, $(\epsilon_2+1-y)(\epsilon_1+y)\ge a_{12}a_{21}$ for $y\in [0,1]$, thus
one can choose a constant $\eta$ such that $(\epsilon_2+1-y)^{-1}a_{21}\le \eta\le (\epsilon_1+y)a_{12}^{-1}$ for $y\in [0,1]$. With $v=(1,\eta e)$, we have
$[D(t)-A(t)]v\ge 0,t\in\R$ and $ [D(t)-A(t)]v\not\equiv 0$. Furthermore, for $u_1,u_2>0$ we have
$$M(t)\left[ \begin{array}{c}u_1\\u_2\end{array}\right]\ge\left[ \begin{array}{c} u_1(\be_1-\epsilon_1+1)+u_2a_{12}e^{-(y+1)}\\
u_1a_{21}e^y+u_2( \be_2-\epsilon_2+1)\end{array}\right].
$$ 
Assume also that
 \begin{equation}\label{3.20}
{\rm either}\ \be_i-\epsilon_i+1\ge 0\ {\rm for\ some}\ i\in \{1,2\}\q {\rm or}\q e^2(\epsilon_1-\be_1-1)(\epsilon_2-\be_2-1)<a_{12}a_{21}.
\end{equation}
If either $\be_1-\epsilon_1+1\ge 0$ or 
$\be_2-\epsilon_2+1\ge 0$, one  finds $u\gg 0$ such that  $M(t)u\gg 0, t\in\R$; if
$\be_i-\epsilon_i+1<0$ for $i=1,2$ and
$e^2(\epsilon_1-\be_1-1)(\epsilon_2-\be_2-1)<a_{12}a_{21},$ then $M(t)u\gg 0, t\in\R$
 with $u=(1,u_2)$ and $u_2$ chosen so that
 $$e^2(\epsilon_1-\be_1-1)a_{12}^{-1}<u_2<a_{21}(\epsilon_2-\be_2-1)^{-1}.$$
 From Theorem \ref{thm3.2}, \eqref{3.19},\eqref{3.20} imply that there is a positive $\pi$-periodic solution of \eqref{3.18}.
\end{exmp}

\section{An application to periodic Nicholson systems}
\setcounter{equation}{0}

For a general system \eqref{0}, it is important to establish conditions for
the global attractivity of the positive $\omega$-periodic solution, whose existence was shown in Theorem \ref{thm3.1}. In order to obtain this global asymptotic behavior of solutions, it is clear that  additional constraints depending strongly on the particular shape of the nonlinearities  $h_{ik}$ should be imposed.
 In this section, we analyze
this situation in the case of a Nicholson system \eqref{1.2}  with constant discrete delays all multiple of the period. For simplicity, we only consider one delay in each equation of the system, 
but straightforward changes allow to consider several delays (all multiple of the period)  as in \eqref{1.2}.

 Consider the periodic Nicholson's system
\begin{equation}\label{4.1}
x_i'(t)=-d_i(t)x_i(t)+\sum_{j=1,j\ne i}^n a_{ij}(t)x_j(t)+ \be_i(t)    x_i(t-m_i\omega)e^{-c_i(t)x_i(t-m_i\omega)},\ i=1,\dots,n,
\end{equation}
where $m_i\in\N,\omega>0,d_i(t),a_{ij}(t), \be_i(t), c_i(t)$ are continuous and $\omega$-periodic, with $d_i(t), \be_i(t), c_i(t)$ positive and $a_{ij}(t)$ nonnegative, for all $i,j$.  

We start with an auxiliary lemma regarding the particular nonlinearity $h(x)=xe^{-x}$.

\begin{lem}\label{lem5.1} For any $x\in (0,2)$, define $G_x:[0,\infty)\to \R$ by
$$G_x(y)=\system{&\frac{h(y)-h(x)}{y-x},\q &y\ne x\cr &(1-x)e^{-x},\q &y=x\cr}$$
where $h(x)=xe^{-x}, x\ge 0$. Then,  for each  $m\in (0,1)$ there is
$\de(x):=\max_{y\ge m} |G_x(y)|<e^{-x}.$
\end{lem}

\begin{pf}  Fix $x\in (0,2)$, and consider $G_x$ defined as above. Note that $G_x(x)=h'(x)$. It was shown in \cite{FariaRost} that
$|h(y)-h(z)|<e^{-z}|y-z|$ for all $y>0$ and $z\in (0,2]$. Since $G_x$ is continuous and $G(\infty)=0$, for any  $m\in (0,1)$ there exists
$\de(x):=\max_{y\ge m} |G_x(y)|$. But $\de(x)<e^{-x}$ because $|G_x(y)|<e^{-x}$ for $y\ne x$ and
$|G_x(x)|=|1-x|e^{-x}<e^{-x}$.  \qed\end{pf}

Next, we denote
$$c_i^-:=\min_{t\in [0,\omega]}c_i(t),\q c_i^+:=\max_{t\in [0,\omega]}c_i(t),\q i=1,\dots,n.$$
%
\begin{lem}\label{lem5.2} 
For some   positive vector $v=(v_1,\dots,v_n)\in \R^n$, suppose that
\begin{equation}\label{5.10}
\al_i(v):=\min_{t\in[0,\omega]} \frac{\be_i(t)v_i}{d_i(t)v_i-\sum_{j\ne i} a_{ij}(t)v_j}>1,\ 1\le i\le n,
\end{equation}
 and define
 \begin{equation}\label{5.10b}
 \ga_i(v):=\max_{t\in[0,\omega]} \frac{\be_i(t)v_i}{d_i(t)v_i-\sum_{j\ne i} a_{ij}(t)v_j},\ 1\le i\le n.
 \end{equation}
 A positive $\omega$-periodic solution   $x^*(t)$ of \eqref{4.1} (whose existence is given in Theorem \ref{thm3.2}) satisfies
$$\frac{x_i^*(t)}{v_i}\le \max_{1\le j\le n} \frac{\log \ga_j(v)}{v_jc_j^-},\ \ t\in [0,\omega],\, i=1,\dots,n.$$
\end{lem}

\begin{pf} Observe that \eqref{5.10} implies that $[D(t)-A(t)]v\gg 0$  and
$(M(t)v)_i\ge (\al_i(v)-1)\eta_i>0$ for $t\in[0,\omega]$, where
$\eta_i=\min_{t\in [0,\omega]}(d_i(t)v_i-\sum_{j\ne i}a_{ij}(t)v_j)$. Hence, \eqref{4.1} satisfies (H1)-(H5).
 Incorporating the scaling $\hat x_i(t)=x_i(t)/v_i\, (1\le i\le n)$ in the coefficients, one further supposes that condition \eqref{5.10} holds with $v={\bf 1}$, and replace $ a_{ij}(t)$ by $\hat a_{ij}(t)=v_i^{-1}a_{ij}(t)v_j$ and $c_i(t)$ by $v_ic_i(t)$ in \eqref{4.1}. Below, we drop the hats in the transformed system, not forgetting however  to insert the weights $v_i's$ in the final estimates.

Since $x^*(t)$ is $\omega$-periodic, it satisfies
\begin{equation}\label{5.10'}
(x_i^*)'(t)=-d_i(t)x_i^*(t)+\sum_{j\ne i} a_{ij}(t)x_j^*(t)+ \be_i(t) x_i^*(t)e^{-v_ic_i(t)x_i^*(t)},\ i=1,\dots,n.
\end{equation}
For $t_0\in [0,\omega]$ and $i\in\{1,\dots,n\}$ such that $ \max_{t\in [0,\omega]}|x^*(t)|=x_i^*(t_0)$, we have 
\begin{equation*}
\begin{split}
0&=-d_i(t_0)x_i^*(t_0)+\sum_{j\ne i} a_{ij}(t_0)x_j^*(t_0)+ \be_i(t_0)x_i^*(t_0)e^{-v_ic_i(t_0)x_i^*(t_0)}\\
&\le \Big(d_i(t_0)-\sum_{j\ne i} a_{ij}(t_0)\Big) x_i^*(t_0)\Big [-1+\ga_i(v)e^{-v_ic_i(t_0)x_i^*(t_0)}\Big]
\end{split}
\end{equation*}
 which implies $e^{v_ic_i(t_0)x_i^*(t_0)}\le \ga_i(v)$, thus $x_i^*(t_0)\le \log \ga_i(v)/(v_ic_i^-)$, and the result follows.\qed\end{pf}
 
 \begin{thm}\label{thm4.1}   For \eqref{4.1}, suppose that there is a vector $v=(v_1,\dots,v_n)\gg 0$ such that
 \begin{equation}\label{5.10''}
 \begin{split}
\al_i(v):=&\min_{t\in[0,\omega]}\frac{\be_i(t)v_i}{d_i(t)v_i-\sum_{j\ne i} a_{ij}(t)v_j}>1\\
 \ga_i(v):=& \max_{t\in[0,\omega]}\frac{\be_i(t)v_i}{d_i(t)v_i-\sum_{j\ne i} a_{ij}(t)v_j}< e^{\frac{2c_0(v)}{c^0(v)}},\q1\le i\le n,
\end{split}
\end{equation}
where $c_0(v)=\min_{1\le i\le n}(v_ic_i^-),\, c^0(v)=\max_{1\le i\le n}(v_ic_i^+)$. Then there exists a unique  positive $\omega$-periodic solution $x^*(t)$, which  is a global attractor of all other positive solutions of \eqref{4.1}; that is,
$x(t)-x^*(t)\to 0$ as $t\to\infty$
for any solution $x(t)=x(t,0,\phi)$ of \eqref{4.1} with initial condition $\phi\in C_0$.
\end{thm}

%

\begin{pf} As before, suppose that \eqref{5.10} holds with $v={\bf 1}$,  and replace $c_i(t)$ by $v_ic_i(t)$, so that  \eqref{4.1} reads as
\begin{equation}\label{5.09}
x_i'(t)=-d_i(t)x_i(t)+\sum_{j=1,j\ne i}^n a_{ij}(t)x_j(t)+ \frac{\be_i(t)}{v_ic_i(t)}    h\big(v_ic_i(t)x_i(t-m_i\omega)\big),\ i=1,\dots,n,
\end{equation}
where $h(x)=xe^{-x}$ as in Lemma \ref{lem5.1}.

We have $\al_i(v)>1,\ \ga_i(v)<e^{2c_0(v)/c^0(v)}$.
From Theorem \ref{thm3.1} and Lemma \ref{lem5.2}, there is a  positive $\omega$-periodic solution $x^*(t)$ of  \eqref{5.09} whose components satisfy
$0<v_ic_i(t)x_i^*(t)\le c^0(v)x_i^*(t)\le c^0(v)\max_i\big (\frac{\log \ga_i(v)}{v_ic_i^-}\big)< 2$ for $t\in [0,\omega]$.  Effecting the change of variables $y_i(t)=\frac{x_i(t)}{x_i^*(t)}-1$ and using \eqref{5.10'},  \eqref{5.09} becomes
\begin{equation}\label{5.11}
\begin{split}
y_i'(t)=\frac{1}{x_i^*(t)}\Big \{&-d_i^*(t)y_i(t)+\sum_{j\ne i} a_{ij}(t)x_j^*(t)y_j(t)\\
&+\frac{\be_i(t)}{v_ic_i(t)}\Big [h\Big(v_ic_i(t)x_i^*(t)(1+y_i(t-m_i\omega))\Big)- h\big (v_ic_i(t)x_i^*(t)\big)\Big ]\Big\},
\end{split}
\end{equation}
where
$$d_i^*(t)=\sum_{j\ne i} a_{ij}(t)x_j^*(t)+\be_i(t)x_i^*(t)e^{-v_ic_i(t)x_i^*(t)}.$$

Let $y(t)=(y_1(t),\dots,y_n(t))$ be any solution of \eqref{5.11} with initial condition $y_0\ge -1, y(0)>-1$. 
Define $-z_i=\liminf_{t\to\infty}y(t), u_i=\limsup_{t\to\infty}y(t)$, and $u=\max_i u_i,z=\max_iz_i$. From the uniform persistence of \eqref{4.1}, $x_i^*(t)(1+y_i(t))\ge m,t\ge 0,$ for some $m\in (0,1)$, and $-1<-z_i\le u_i<\infty$. 

It is sufficient to show that $\max(u,z)=0$. Suppose that $\max(u,z)=u>0$ (the situation  $\max(u,z)=z$ is treated in a similar way).  Choose $i$ such that $u=u_i$ and take a sequence $t_k\to\infty$ with $y_i(t_k)\to u, y_i'(t_k)\to 0$.  Let $\vare>0$ be small. From \eqref{5.11} and Lemma \ref{lem5.1}, for large $k$ we get 
\begin{equation}\label{5.12}
\begin{split}
y_i'(t_k)&\le \frac{1}{x_i^*(t_k)}\bigg [-d_i^*(t_k)+\sum_{j\ne i} a_{ij}(t_k)x_j^*(t_k)\Big ]y_i(t_k)\\
&\hskip 3.2cm+\be_i(t_k)x_i^*(t_k)\de \Big(v_ic_i(t_k)x_i^*(t_k)\Big)  |y_i(t_k-m_i\omega)|\bigg]+O(\vare)\\
&=\be_i(t_k) \Big [ -e^{-v_ic_i(t_k)x_i^*(t_k)}y_i(t_k)+\de \Big(v_ic_i(t_k)x_i^*(t_k)\Big)  |y_i(t_k-m_i\omega)|\Big ] +O(\vare).
\end{split}
\end{equation}

For some subsequence of $(t_k)$, still denoted by $(t_k)$,  $\lim_k v_ic_i(t_k)x_i^*(t_k)=\xi\in (0,2),\lim_k \be_i(t_k)=b>0, \lim_ky_i(t_k-m\omega)=w\in [-z,u]$. The estimate \eqref{5.12} leads to
$$0\le b(-e^{-\xi}u+\de(\xi)|w|)\le b(-e^{-\xi}+\de(\xi))u$$
which is not possible because $\de(\xi)<e^{-\xi}$ for any $\xi\in (0,2)$. Thus $u=0$.\qed\end{pf}

Several important consequences can  be deduced from Theorem \ref{thm4.1}. 

\begin{cor}\label{cor4.1} 
 Consider the classic periodic Nicholson's equation with a delay multiple of the period:
\begin{equation}\label{5.16}
x'(t)=-d(t)x(t)+\be(t)x(t-m\omega)e^{-c(t)x(t-m\omega)},
\end{equation}
where $\omega>0,d(t),\be(t), c(t)$ are continuous, positive and $\omega$-periodic functions and $m\in\N$.  
Set
$\min_{t\in [0,\omega]}c(t)=c^-,\ \max_{t\in [0,\omega]}c(t)=c^+,$ and
suppose that 
$$1<\frac{\be(t)}{d(t)}<e^{2c^-/c^+},\q t\in[0,\omega].$$
Then there exists a unique  positive $\omega$-periodic solution $x^*(t)$, which  is a global attractor of all other positive solutions of \eqref{5.16}. In particular, if $c(t)\equiv c>0$, the global attractivity of $x^*(t)$ holds true if $1<\be(t)/d(t)<e^2,t\in[0,\omega].$
\end{cor}

\begin{cor}\label{cor4.2}   Consider the  autonomous Nicholson's system
\begin{equation}\label{5.17}
x_i'(t)=-d_ix_i(t)+\sum_{j=1,j\ne i}^n a_{ij}x_j(t)+\sum_{k=1}^m \be_{ik}    x_i(t-\tau_{ik})e^{-c_ix_i(t-\tau_{ik})},\
i=1,\dots,n,\ t\ge 0,\end{equation}
where $d_i>0,  c_i>0, a_{ij}\ge 0\, (j\ne i), \tau_{ik}\ge 0, \be_{ik}\ge 0$ with $\be_i:=\sum_{k=1}^m \be_{ik} >0$ for all $i,j,k$.

If there exists a vector $v\gg 0$ such that
\begin{equation}\label{5.17'}
1<\gamma_i(v)<e^{\frac{2\min_j (v_jc_j)}{\max_j(v_jc_j)}}\q {\rm for}\q \gamma_i(v)=\frac{\be_iv_i}{d_iv_i-\sum_{j\ne i}a_{ij}v_j}\, ,\q i=1,\dots,n,
\end{equation}
then there exists a unique  positive equilibrium which  is a global attractor of all  positive solutions of \eqref{5.17}.\end{cor}

\begin{pf}
The proof follows as the proof of Theorem \ref{thm4.1}, with the positive $\omega$-periodic solution $x^*(t)$ replaced by the positive equilibrium $x^*$.\qed\end{pf}

\begin{cor}\label{cor4.3}   Consider the  autonomous Nicholson's system
\begin{equation}\label{5.18}
x_i'(t)=-d_ix_i(t)+\sum_{j=1,j\ne i}^n a_{ij}x_j(t)+\sum_{k=1}^m \be_{ik}    x_i(t-\tau_{ik})e^{-x_i(t-\tau_{ik})},\
i=1,\dots,n,\ t\ge 0,\end{equation}
where $d_i>0,  a_{ij}\ge 0\, (j\ne i), \tau_{ik}\ge 0, \be_{ik}\ge 0$ with $\be_i:=\sum_{k=1}^m \be_{ik} >0$ for all $i,j,k$.
If there exists a vector $v\gg 0$ such that
\begin{equation}\label{5.19}
1<\gamma_i(v)<e^{\frac{2\min_j (v_j)}{\max_j(v_j)}},\q i=1,\dots,n,
\end{equation}
where $\gamma_i(v)$ are defined in \eqref{5.17'}, 
then there exists a unique  positive equilibrium which  is a global attractor of all  positive solutions of \eqref{5.18}. In particular, this is the case if
\begin{equation}\label{5.20}
1<\ga_i<e^2\q {\rm for}\q \ga_i: =\frac{\be_i}{d_i-\sum_{j\ne i}a_{ij}},\q  i=1,\dots,n.
\end{equation}
\end{cor}

\begin{rmk}\label{rmk5.3} For the particular situation \eqref{5.18}, the result in Corollary \ref{cor4.3}  was proven in \cite{FariaRost} under the  hypothesis $1<\ga_i\le e^2,\,1\le i\le n$, for $\gamma_i=\gamma_i({\bf 1})$ as in \eqref{5.20}. To obtain the result for $\max_i\ga_i=e^2$, the proof however uses results on $\omega$-limit sets for {\it autonomous} DDEs, which do not carry out for \eqref{4.1}, much less for more general periodic systems  \eqref{1.2}. On the other hand, adapting the proof in \cite{FariaRost}, it is now apparent that Corollary \ref{cor4.2}  is valid with
$1<\frac{\be_iv_i}{d_iv_i-\sum_{j\ne i}a_{ij}v_j}\le \exp\big(\frac{2\min_j (v_jc_j))}{\max_j(v_jc_j)}\big),\ i=1,\dots,n,$ which improves the criterion in \cite{FariaRost}. 
\end{rmk}

\begin{exmp} Consider the 2-dimensional $\omega$-periodic Nicholson system
  with one single discrete delay given by
\begin{equation}\label{1.9}
\begin{split}
&x_1'(t)=-a_1(t)x_1(t)+b_1(t)x_2(t)+c_1(t)x_1(t-\tau)e^{-x_1(t-\tau)}\\
&x_2'(t)=-a_2(t)x_2(t)+b_2(t)x_2(t)+c_2(t)x_2(t-\tau)e^{-x_2(t-\tau)},
  \end{split}
  \end{equation}
 where $a_i(t),b_i(t), c_i(t)\, (i=1,2)$ are positive, continuous and $\omega$-periodic functions, and  suppose that $\tau=m\omega$ for some $m\in\N$. Applying  Theorem \ref{thm4.1}, we derive that  \eqref{1.9} has a globally attractive positive $\omega$-periodic solution if there exist positive constants $v_1,v_2$ such that
$$1<  \frac{c_1(t)v_1}{a_1(t)v_1-b_1(t)v_2}<e^2\ {\rm and} \ 1<\frac{c_2(t)v_2}{a_2(t)v_2-b_2(t)v_1}<e^2,\q t\in[0,\omega].$$
In particular, this assertion is valid if
\begin{equation}\label{4.15}
1<\frac{c_i(t)}{a_i(t)-b_i(t)}<e^2,\q t\in[0,\omega],i=1,2.
\end{equation}
A similar result holds with  $\tau,\omega$  rationally dependent. 

We now compare this criterion with the one in   \cite{Troib}.  Recently, Troib   \cite{Troib}  used the continuation theorem of coincidence degree to show the  existence of a positive periodic solution $x^*(t)$   for \eqref{1.9} under the following constraints:
 \begin{equation}\label{1.9'}2D_i\min \{ e^{A_1}, e^{A_2}\}<A_i\le 4D_i \max \{ e^{A_1}, e^{A_2}\},\ 2C_i>e^{A_i}A_i,\q i=1,2,
 \end{equation}
where $A_i=2 \omega \ol a_i, B_i=\omega \ol b_i, C_i=\omega \ol c_i, D_i=\max \{ B_i,C_i\}, \, i=1,2$, and the notation $\ol f=\frac{1}{\omega}\int_0^\omega f(t)\,dt$ is used for an $\omega$-periodic function.  We observe however  that for the particular case of the scalar periodic Nicholson equation \eqref{NP}  with $\tau(t)\equiv \tau$ and $c(t)\equiv 1$, the criterion in  \cite{Troib} does not apply, since the  conditions \eqref{1.9'}  would read as $\ol \be >\ol d e^{2\omega \ol d},\ol \be e^{2\omega \ol d}<\ol d<2 \ol \be e^{2\omega \ol d}$, and the set of functions $\be, d$ satisfying these conditions is empty. 
By using a suitable Lyapunov functional, in   \cite{Troib} the author
  further obtained the global asymptotic stability of $x^*(t)$ under  the additional restrictions
 \begin{equation*}
\min_{t\in\R} c_i(t)>e^{M_i} \max_{t\in\R}a_i(t),
\max_{t\in\R}c_i(t)<\big(\min_{t\in\R} a_i(t)-\max_{t\in\R}b_i(t)\big)e^2,\ i=1,2,
  \end{equation*}
for $M_i\ge\limsup_{t\to\infty}x_i(t),\, i=1,2$, for all  solutions $(x_1(t),x_2(t))$ of \eqref{1.9},
a hypothesis much stronger than the assumptions \eqref{4.15}.
\end{exmp}

\section*{Acknowledgement} 
This work was partially supported by Funda\c c\~ao para a Ci\^encia e a Tecnologia under project UID/MAT/04561/2013.


\section*{References} 
 
 \begin{thebibliography}{99}

{


\bibitem{AmsterIdels} P. Amster, L. Idels, Periodic solutions in general scalar non-autonomous models with delays, Nonlinear  Differential Equations Appl. 20 (2013), 1577--1596.





%
  \bibitem{BIT}
   L. Berezansky, L. Idels, L. Troib,  Global dynamics of Nicholson-type delay systems with applications,
Nonlinear Anal. Real Word Appl. 12 (2011), 436--445.

\bibitem{Chen}
Y. Chen, Periodic solutions of delayed periodic Nicholson's blowflies models, Can. Appl. Math. Q.  11 (2003), 23--28.




\bibitem{DingNieto} H.-S. Ding,  J.J. Nieto, A new approach for positive almost periodic solutions to a class of Nicholson's blowflies model,  J. Comput. Appl. Math. 253 (2013), 249--254.



\bibitem{Diekmann01}
O. Diekmann, M. Gyllenberg, H. Huang, M. Kirkilionis,
J.A.J. Metz, H.R. Thieme, On the formulation and analysis of general
deterministic structured population models
II. Nonlinear theory, J. Math. Biol. 43  (2001), 157--189.
%
%
%
%


\bibitem{Faria11} T. Faria, Global asymptotic behaviour for a Nicholson model
 with patch structure and multiple delays, Nonlinear Anal. 74 (2011), 7033--7046.




 
 \bibitem{FOS}
T. Faria, R. Obaya, A.M. Sanz, Asymptotic behaviour for non-monotone delayed perturbations of monotone non-autonomous linear ODEs, {\it submitted} (2016). 



\bibitem{FariaRost}
T. Faria,  G. R\"ost, Persistence, permanence and global stability
or an $n$-dimen\-sional Nicholson system,  J. Dynam. Differential Equations 26 (2014), 723--744.



\bibitem{Fiedler}   
M. Fiedler,  Special Matrices and Their Applications in Numerical Mathematics, Martinus Nijhoff Publ. (Kluwer), Dordrechit, 1986.

\bibitem{FLT} D. Franco, E. Liz, P.J. Torres, Existence of periodic solutions for functional equations with periodic delay, Indian J. Pure Appl. Math. 38 (2007), 143--152.

\bibitem{GBN}
W.S.C.  Gurney,  S.P. Blythe, R.M. Nisbet,   Nicholson's blowflies revisited,  Nature  287 (1980), 17--21.

\bibitem{Hale}
J. K.  Hale,
Asymptotic Behavior of Dissipative Systems, Math. Surveys Monogr., Vol. 25, Amer. Math. Soc., Providence, RI, 1988.

\bibitem{HaleLunel}
J.K. Hale, S.M.  Verduyn Lunel,  Introduction to Functional Differential
Equations, Springer-Verlag, New-York, 1993.

\bibitem{Hofbauer} 
 J. Hofbauer, An index theorem for dissipative systems, Rocky Mountain J. Math.  20 (1990), 1017--1031.
 




\bibitem{Kuang}  
 Y. Kuang, Delay Differential Equations with Applications in Population dynamics, Academic Press, New York, 1993. 
 



\bibitem{LiDu}
 J. Li and C. Du, Existence of positive periodic solutions for a generalized Nicholson's blowflies model, J. Comput. Appl. Math. 221 (2008), 226--233.
 
 \bibitem{Liu}
 B. Liu, Global stability of a class of delay differential systems, J. Comput. Appl. Math. 233 (2009), 217--223.


 
 \bibitem{Liu11}
 B. Liu, The existence and uniqueness of positive periodic solutions of Nicholson-type delay systems, Nonlinear Anal. Real Word Appl. 12 (2011), 3145--3451.

\bibitem{Liu14} 
B. Liu, Global exponential stability of positive periodic solutions for a delayed Nicholson's blowflies model, J. Math. Anal. Appl. 412 (2014), 212--221.

\bibitem{MetzDiek}
J.A.J. Metz, O. Diekmann, The Dynamics of Physiologically Structured Populations (Eds.), Lecture Notes in Biomath.  68, Springer-Verlag, 1986.


\bibitem{Nicholson} A.J. Nicholson, An outline of the dynamics of animal populations, Austral. J. Zool. 2 (1954), 9--65.


\bibitem{ObayaSanz}
R. Obaya, A.M. Sanz, Uniform and strict persistence in monotone skew-product semiflows with applications to non-autonomous Nicholson systems, J. Differential Equations 261 (2016),  4135--4163.



\bibitem{SakerAg}
S.H. Saker and S. Agarwal, Oscillation and global attractivity in a periodic
Nicholson's blowflies model, Math. Comput. Modelling 35 (2002), 719--731.

\bibitem{Smith}   H.L. Smith,   Monotone Dynamical Systems. An
Introduction to the Theory of Competitive and Cooperative Systems,
 Mathematical Surveys and Monographs, Amer.~Math.~Soc., Providence, RI, 1995.
 
 \bibitem{Smith11}  H.L. Smith,   An Introduction to Delay Differential Equations with Applications to Life Sciences, Texts in Applied Mathematics Vol. 57, Springer, Berlin, 2011.
 


 



\bibitem{Troib}
L. Troib, Periodic solutions of Nicholson-type delay differential systems, Funct. Diff. Equations 21 (2014),  171--187.

\bibitem{Wang04} H. Wang, Positive periodic solutions of functional differential equations, J. Differential Equations 202 (2004), 354--366.

\bibitem{Wang13}  L. Wang, Almost periodic solution for Nicholson's blowflies model with patch structure and linear harvesting terms, Appl. Math. Model. 37 (2013),  2153--2165.


\bibitem{WWC}
W. Wang, L. Wang, W. Chen, Existence and exponential stability of positive almost periodic solution for Nicholson-type delay systems,
Nonlinear Anal. Real Word Appl. 12 (2011), 1938--1949.

%



\bibitem{ZDC}
N. Zhang, B. Dai, Y. Chen, Positive periodic solutions of nonautonomous functional differential systems,
 J. Math. Anal. Appl. 333 (2007), 667--678.

\bibitem{Zhou}
 Q. Zhou, The positive periodic solution for Nicholson-type delay system with linear harvesting terms, Appl. Math. Model. 37 (2013), 5581--5590.

}
\end{thebibliography}
\end{document}